\documentclass[12pt,a4paper]{amsart}

\usepackage{amsaddr}
\usepackage{amssymb}
\usepackage[final]{showkeys}
\usepackage{microtype}
\usepackage{color}
\usepackage{graphicx}
\usepackage[hmargin=3cm,vmargin={3.5cm,4cm}]{geometry}

\numberwithin{equation}{section}

\allowdisplaybreaks

\begin{document}

    \title[G-fractional diffusion on bounded domains]{G-fractional diffusion  on bounded domains in $\mathbb{R}^d$}

\author{L. Angelani$^{1,2}$ and R. Garra$^{3}$}
\address{$^1$ Istituto dei Sistemi Complessi, Consiglio Nazionale delle Ricerche, P.le A. Moro 2, 00185 Roma, Italy}
\address{$^2$ Dipartimento di Fisica, Sapienza Universit\`a di Roma, P.le A. Moro 2, 00185 Roma, Italy}
\address{$^3$ Section of Mathematics, International Telematic University Uninettuno, Corso Vittorio Emanuele II, 39, 00186 Roma, Italy}

    \date{\today}

    \begin{abstract}
     In this paper we study $g$-fractional diffusion on bounded domains in $\mathbb{R}^d$ with absorbing boundary conditions. A new general and explicit representation of the solution is obtained. We study the first-passage time distribution, showing the dependence on the particular choice of the function $g$. Then, we specialize the analysis to the interesting case of a rectangular domain. Finally, we briefly discuss the connection of this general theory with the physical application to the so-called fractional Dodson diffusion model, recently discussed in the literature. 
    
    \end{abstract}

       \keywords{Fractional diffusion equation, first-passage time, $g$-fractional diffusion in bounded domain}

    \maketitle
	
	\section{Introduction}

  Time-fractional diffusion processes are widely studied in the literature for their 
  relation to continuous-time random walks (see, e.g., \cite{ksbook}) and for pervasive applications in various fields, from applied physics to pure mathematics and probability (see, for example, the recent monograph \cite{lenzi} and references therein). More recently, $g$-fractional diffusive equations based on the application of fractional derivatives with respect to another function (also named, in the mathematical literature, $\Psi$-fractional derivatives) 
  have gained greater interest. In particular, in a series of interesting papers
  \cite{kosto, kosto1, kosto2}, the authors discussed the relevance of this approach for physical models of anomalous diffusion. 
  In the mathematical literature, starting from the paper by Almeida \cite{alm},
  many papers have been devoted to the analysis of fractional differential equations with respect to another function and this is a developing field, 
  as it allows nontrivial generalizations of classical equations involving Caputo derivatives.\\

  In \cite{PRE2023} we have considered the one-dimensional $g$-fractional diffusion equation with absorbing boundary conditions. An interesting outcome of the analysis developed in \cite{PRE2023} is that the explicit solution can be found and particular choices of the $g$-function leads to a finite mean first-passage time (MFPT), differently from the time-fractional diffusion in bounded domain involving the classical Caputo derivative (which is the special case $g(t) = t$). In this paper we consider, for the first time, the $g$-fractional diffusion in $d$-dimensional bounded domains with absorbing boundary conditions. We obtain the explicit representation of the solution for the Dirichlet problem in $\mathbb{R}^d$, which can be applied to several particular diffusive models in bounded domains. 
  We emphasize that the main difference with the more particular (and simple) one-dimensional case previously treated in \cite{PRE2023} lies in the broad generality of the main results presented here, which can be used for realistic diffusive models in higher dimensions.
  We then study the first-passage time distribution (FPTD), discussing the 
  condition for a finite MFPT.\\
 As a first application of the general results obtained here, we consider the special and interesting case of $g$-fractional diffusion in rectangular domains. We recall that fractional diffusions in multidimensional rectangular
  domains have been the subject of recent interest in the mathematical literature
  (see, for example, \cite{psu}). \\
  In the final section, 
  we investigate the case of 
  fractional Dodson equation, recently introduced in \cite{garra}. This is an interesting heuristic model of fractional diffusion that includes, at the same time, memory effects and the particular exponential time-dependence of the diffusion coefficient. The peculiarity of this model relies on the fact that the 
  corresponding $g$-function is upper bounded, resulting in the existence of stationary solutions, 
  finite values of the asymptotic survival probability and, therefore, undefined MFPT.

    \section{G-fractional diffusion in bounded domains}

We consider a $g$-fractional diffusion process on a bounded domain $\Omega$ in $\mathbb{R}^d$:
\begin{equation}
\left(\frac{^C\partial_g^\alpha u}{\partial t^{\alpha}}   \right) ({\bf r},t)
= D \nabla^2 u({\bf r},t),
\label{GFDE}
\end{equation}
where ${\bf r}=\{x_1, ..., x_d\}$,
$\nabla^2 = \sum_{i=1}^d \frac{\partial^2}{\partial x_i^2}$ is the Laplacian
in dimension $d$, $D$ is the generalized diffusion constant and the 
Caputo-type $g$-fractional derivative of order $\alpha \in (0,1)$ is defined as
    \begin{equation}
    \left(\frac{^C\partial_g^\alpha u}{\partial t^{\alpha}}   \right) ({\bf r},t)
    := \frac{1}{\Gamma(1-\alpha)} \int_0^t (g(t)-g(\tau))^{-\alpha}\ 
    \frac{\partial u}{\partial \tau} ({\bf r},\tau) d\tau,
    \end{equation} 
where $g(t)$ is a deterministic function such that $g(0) = 0 $ and $g'(t)>0$ for $t>0$, where we denote by $g' = dg/dt$ the first order time derivative. 
We consider generic initial condition
\begin{equation}
    u({\bf r},0) = u_0({\bf r}) ,
\end{equation}
and absorbing (Dirichlet) boundary conditions 
\begin{equation}
    u({\bf r},t) = 0 ,\hspace{1cm} {\bf r} \in \partial \Omega.
\label{bc}
\end{equation}
Solutions can be found  by the method of separation of variables, i.e.
$$u({\bf r},t) = X({\bf r}) T(t).$$
We have to solve the two equations 
\begin{equation}
\label{Teq}
\frac{^C d_g^\alpha T(t)}{d t^{\alpha}}
 = -\lambda D T(t),
 \end{equation}
and 
\begin{equation}
 \label{Xeq}
\nabla^2 X({\bf r}) = -\lambda X({\bf r}).
\end{equation}
The solution of (\ref{Teq}) is 
\begin{equation}
T(t) = T(0) \ E_\alpha(-\lambda D g(t)^\alpha),
\end{equation}
where $E_\alpha(\cdot)$ denotes the one-parameter Mittag-Leffler function \cite{ml}
\begin{equation}
E_{\alpha}(x) = \sum_{k=0}^\infty \frac{x^k}{\Gamma(\alpha k+1)}.
\end{equation}
The eigenvalue problem (\ref{Xeq}) -- with boundary conditions (\ref{bc}) -- is solved by an infinite sequence of pairs $(\lambda_n,\phi_n)$, with $n\geq 1$ ($\lambda_1<\lambda_2 < ...$) and  
$\phi_n({\bf r})$ is a sequence of functions that form a complete orthonormal set in $L^2(\Omega)$ \cite{Gre2013, Nane1, Nane}.\\
The solution of the $g$-fractional diffusive equation (\ref{GFDE})
can then be expressed as 
\begin{equation}
u({\bf r},t) = \sum_{n=1}^\infty u_{0,n} \ \phi_n({\bf r}) \ 
E_\alpha(-\lambda_n D g(t)^\alpha),
\label{urt}
\end{equation}
where 
\begin{equation}
\label{u0n}
u_{0,n} = \int_{\Omega} d{\bf r} \ \phi_n({\bf r})\ u_0({\bf r}).
\end{equation}

\section{first-passage times}
We study here first-passage problems. 
The FPTD $\varphi(t)$ is defined by
\begin{equation}
    \varphi(t) = - \frac{d {\mathbb P}}{dt}(t),
\label{FPTD}
\end{equation}
where ${\mathbb P}(t)$ is the survival probability, i.e. the probability that a particle has not been absorbed until time $t$
\begin{equation}
    {\mathbb P}(t) = \int_{\Omega} d{\bf r} \ u({\bf r},t).
    \label{SP}
\end{equation}
We are assuming here that the survival probability goes to zero for $t \to \infty$, i.e., 
the particle will surely be absorbed during the entire process.
This corresponds to consider functions $g(t)$ such that $\lim_{t \to \infty} g(t) = +\infty$.
In the last section, treating the fractional Dodson diffusion, we will discuss the main 
consequences of relaxing such an assumption.\\
The MFPT $\tau$ is the first moment of (\ref{FPTD})
\begin{equation}
    \tau = \int_0^\infty dt \ t \ \varphi(t).
    \label{MFPT}
\end{equation}
By using (\ref{urt}), (\ref{FPTD}) and (\ref{SP}), we can then express the FPTD for $g$-fractional diffusion processes as
\begin{equation}
\varphi(t) = - \sum_{n=1}^\infty u_{0,n} \ \Phi_n \ 
\frac{d}{dt} E_\alpha(-\lambda_n D g(t)^\alpha),
\end{equation}
where 
\begin{equation}
\Phi_n = \int_{\Omega} d{\bf r} \ \phi_n({\bf r}).
\label{Phi}
\end{equation}
By using the property \cite{ml}
\begin{equation}
    E_{\alpha,\alpha}(-x) = - \alpha \frac{d}{dx} E_\alpha(-x),
\end{equation}
where we have introduced the two-parameters Mittag-Leffler function 
\begin{equation}
E_{\alpha,\beta} (x) = \sum_{k=0}^\infty \frac{x^k}{\Gamma(\alpha k + \beta)} ,
\end{equation}
we finally arrive at the expression for the FPTD
\begin{equation}
\varphi(t) = D g'(t) g(t)^{\alpha-1} \sum_{n=1}^\infty \lambda_n u_{0,n} \ \Phi_n \ 
E_{\alpha,\alpha} (-\lambda_n D g(t)^\alpha) .
\label{FPTD1}
\end{equation}
It is interesting to show the long time behavior of $\varphi(t)$.
By using the asymptotic expansion of the Mittag-Leffler function for $|z| \to \infty$ and $\Re(z)<0$
(see \cite{ml}, p. 75)
\begin{equation}
E_{\alpha, \alpha}(z)= - \frac{z^{-2}}{\Gamma(-\alpha)} 
+ O(|z|^{-3}),
\end{equation}
we have that the asymptotic behavior of (\ref{FPTD1}) is
\begin{equation}
\varphi(t) \sim -\frac{g'(t) g(t)^{-(\alpha+1)}}{D \Gamma(-\alpha)} \sum_{n=1}^\infty \frac{u_{0,n} \ \Phi_n}{\lambda_n} ,
\hspace{1.3cm} t \to \infty.
\label{FPTDasy}
\end{equation}
It is worth noting that the above asymptotic form allows the MFPT (\ref{MFPT}) to be finite only for those functions $g(t)$ 
satisfying $\lim_{t\to\infty} t^2 g' g^{-\alpha-1} = 0$, regardless of boundary shape. In other words, finite MFPTs are obtained if $g(t)$ grows asymptotically faster than $t^{1/\alpha}$, similarly to the one-dimensional case \cite{PRE2023}.

\section{Rectangular domains}
We specialize here to the case of rectangle-like domains in
${\mathbb R}^d$, $\Omega=[0,L_1] \times \dots \times [0,L_d]$.
Variable separation allows the eigenfunctions to be written as (by using the multiple index $n=\{n_1,...,n_d\}$)
\begin{equation}
\phi_n({\bf r}) =   \prod_{i=1}^{d}
\phi_{n_i}^{(i)}(x_i) ,
\end{equation}
and the eigenvalues as
\begin{equation}
\lambda_n = \sum_{i=1}^d \lambda_{n_i}^{(i)}.
\label{eigenval}
\end{equation}
Considering absorbing boundary conditions (\ref{bc}) we have 
\cite{Gre2013, Nane}
\begin{eqnarray}
\phi_{n_i}^{(i)}(x) &=& (2/L_i)^{1/2} \ \sin{(\pi n_i x/L_i)} ,\\
\lambda_{n_i}^{(i)} &=& \pi^2 n_i^2 / L_i^2 ,
\label{lamb}
\end{eqnarray}
with $i=1,...,d$ and $n_i\geq1$.
Inserting in (\ref{urt}) we obtain the solution of the $g$-fractional diffusion  equation in rectangular domains with generic initial conditions
\begin{equation}
u({\bf r},t) = 
\sum_{n=1}^\infty u_{0,n}
\left[
\prod_{i=1}^d \sqrt{\frac{2}{L_i}}
\sin{\left(\frac{\pi n_i x_{i}}{L_i}\right)}
\right]
\ E_{\alpha} \left(
-\lambda_n D g(t)^\alpha
\right),
\label{urt1}
\end{equation}
where $\lambda_n$ are given by (\ref{eigenval}).
In the following we consider $\delta$-peaked initial conditions
\begin{equation}
u_0({\bf r}) = \delta({\bf r}-{\bf r}_0),
\end{equation}
where ${\bf r}_0 = \{x_{1,0},...,x_{d,0}\}$. 
We have, 
from (\ref{u0n}),
\begin{equation}
u_{0,n} = \prod_{i=1}^d \phi_{n_i}^{(i)}(x_{i,0}) = 
\prod_{i=1}^d \sqrt{\frac{2}{L_i}} \sin{\left(\frac{\pi n_i x_{i,0}}{L_i}\right)},
\end{equation}
and the solution (\ref{urt1}) reads
\begin{equation}
u({\bf r},t) = 
\sum_{n=1}^\infty 
\left[
\prod_{i=1}^d \frac{2}{L_i}
\sin{\left(\frac{\pi n_i x_{i,0}}{L_i}\right)}
\sin{\left(\frac{\pi n_i x_{i}}{L_i}\right)}
\right]
\ E_{\alpha} \left(
-\lambda_n D g(t)^\alpha
\right).
\label{urt2}
\end{equation}
We now turn to analyze the FPTD (\ref{FPTD1}).
We first note that the terms $\Phi_n$ (\ref{Phi}) can be written as
\begin{equation}
\Phi_n =   \prod_{i=1}^{d} \Phi_{n_i}^{(i)} ,
\end{equation}
where 
\begin{eqnarray}
\Phi_{n_i}^{(i)} &=&   \int_0^{L_i} dx \ \phi_{n_i}^{(i)}(x) =
\sqrt{\frac{2}{L_i}} \int_0^{L_i} dx \ \sin{(\pi n_i x /L_i)} \nonumber\\
&=&
\frac{2 \sqrt{2L_i}}{\pi n_i} , \quad \quad \mbox{if $n_i$ is odd},
\end{eqnarray}
and null otherwise.
We can finally express the FPTD as
\begin{equation}
\varphi(t) = \frac{2^{2d}D g'(t) g(t)^{\alpha-1}}{\pi^d} 
\sum_{n=0}^\infty \lambda_{2n+1} 
\left(
\prod_{i=1}^d \frac{\sin{(\pi (2n_i+1) x_{i,0}/L_i)}}{2n_i+1}
\right)
\ E_{\alpha,\alpha} (-\lambda_{2n+1} D g(t)^\alpha),
\label{FPTD2}
\end{equation}
where
$\lambda_{2n+1} = \sum_{i=1}^d 
\pi^2 (2 n_i+1)^2/L_i^2$.
We note that for $d=1$ the above expression reduces to that obtained in \cite{PRE2023} for the one-dimensional case. 

\begin{figure}[t!]
\includegraphics[width=.78\linewidth] {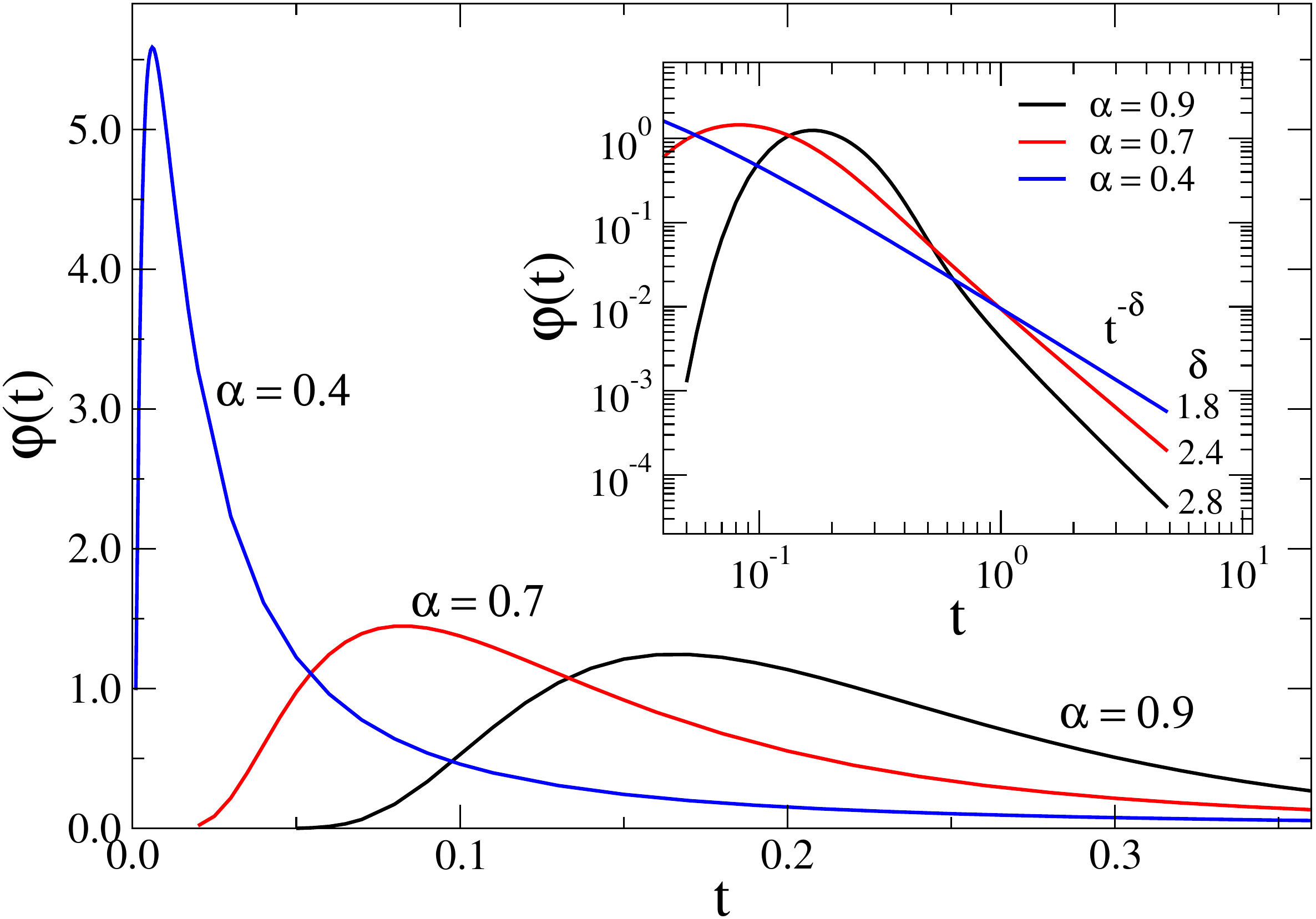}
\caption{
Example of first-passage time distributions at different fractional derivative order $\alpha$ for a particular choice of the $g$-function, $g(t)=t^2$ (Erd\'erlyi-Kober derivative). Inset: the same as in the main panel in logarithmic scale, in order to highlight the long time behavior $t^{-\delta}$.
Only the cases for which $\delta=1+2\alpha>2$ correspond to finite values of the MFPT
(see (\ref{FPTDasy}) and subsequent discussion).
The curves are obtained by numerical evaluation of the expression (\ref{FPTD2}), considering a square box ($d=2$) of size $L=1$ and setting $D=1$, $x_{1,0}=x_{2,0}=L/2$, and $g(t)=t^2$.
}
\label{fig1}
\end{figure}

As an example, Figure \ref{fig1} shows the typical behavior of FPTDs at different values of the fractional order $\alpha$. The curves shown are obtained by numerical evaluation of (\ref{FPTD2}) and correspond to the Erd\'erlyi-Kober derivative
$g(t)=t^\beta$, with $\beta=2$, in a two-dimensional square box with absorbing boundaries.
It is evident the power-law decay at long time, affecting the existence of 
finite MFPTs -- see (\ref{FPTDasy}) and the discussion at the end of the previous section.

\section{The fractional Dodson diffusion}

The Dodson diffusion equation arises in the context of cooling processes in geology \cite{dodson} and 
     it takes the form (see e.g. \cite{crank}, pag.104-105)
     \begin{equation}\label{0}
     \frac{\partial u}{\partial t} = D_0 \, \exp(-\beta t) \, \frac{\partial^2 u}{\partial x^2}, 
     \end{equation}
     where $1/\beta$ is the so called \textit{relaxation time}. 
     This is an interesting model where the diffusivity coefficient is time-dependent and, more precisely, it is an exponentially decreasing function of time.
     As we will see, this slowing down of dynamics generates finite stationary solutions, with important consequences on first-passage  processes.
     
     In the recent paper \cite{garra}, a new generalization of the Dodson diffusion equation was suggested in view of the relevance of the fractional approach for diffusive models with memory effects. 
     In this paper, the fractional Dodson equation is essentially a $g$-fractional diffusion with 
     \begin{equation}
     g(t) =\frac{1 - e^{-\beta t}}{\beta}.
     \label{gDod}
     \end{equation}
    This means that the time-fractional operator appearing in the evolution equation is given by 
\begin{equation}\label{3.1}
           {}^C \left(e^{\beta t} \, \frac{\partial}{\partial t}\right)^\alpha  u(x,t) =\frac{1}{\Gamma(1-\alpha)}\int_0^t \left(\frac{e^{-\beta \tau} - e^{-\beta t}}{\beta} \right) ^{- \alpha} \, \frac{\partial u}{\partial \tau} \, d \tau,
    \end{equation}
       where we used the notation of \cite{garra} to underline that this 
       $g$-fractional derivative physically corresponds to an operator that includes the time-dependence of the diffusivity and the memory effects.
       Moreover, for $\alpha = 1$ we recover the classical equation.
       This is an exploratory generalization of the Dodson diffusion and a concrete potential application of the $g$-fractional approach to diffusive models. \\
       We note that the $g(t)$ (\ref{gDod}) is a bounded function, as
       $\lim_{t \to \infty} g(t) = 1/\beta$.
       Having this in mind, we can apply some of the results obtained in the previous sections to the present case of fractional Dodson diffusion in bounded domains with absorbing boundary conditions. 
        We recall that in \cite{garra}, the authors derived the fundamental solution of the fractional Dodson equation in the 
       one-dimensional case, while the diffusive problem in a bounded domain
       and in higher dimensions has not been considered before in the literature.\\
       Let us consider the $d$-dimensional fractional Dodson equation 
       
\begin{equation}\label{neww}
 	     {}^C \left(e^{\beta t}\frac{\partial}{\partial t}\right)^\alpha u({\bf r},t) =  D \nabla^2 u({\bf r},t) ,
 	     \end{equation}
       under the initial condition
\begin{equation}
    u({\bf r},0) = u_0({\bf r}) ,
\end{equation}
and absorbing boundary conditions
\begin{equation}
    u({\bf r},t) = 0 ,\hspace{1cm} {\bf r} \in \partial \Omega.
\label{bcdD}
\end{equation}
Then, using the results obtained in the previous sections, we have that the solution 
can be expressed as 
\begin{equation}
u({\bf r},t) = \sum_{n=1}^\infty u_{0,n} \ \phi_n({\bf r}) \ 
E_\alpha\bigg(-\lambda_n D \left(\frac{1 - e^{-\beta t}}{\beta}\right)^\alpha\bigg),
\label{urtD}
\end{equation}
where 
\begin{equation}
\label{u0nD}
u_{0,n} = \int_{\Omega} d{\bf r} \ \phi_n({\bf r})\ u_0({\bf r}).
\end{equation}
Thus, the upper bounded $g(t)$ implies the existence of a stationary solution $u_{st.}$:
\begin{equation}
    u_{st.}({\bf r}) = \lim_{t \to \infty} u({\bf r},t) = 
    \sum_{n=1}^\infty u_{0,n} \ \phi_n({\bf r}) \ 
E_\alpha(-\lambda_n D \beta^{-\alpha}).
\label{urtDst}
\end{equation}
We then conclude that the survival probability has a finite asymptotic value
\begin{equation}
{\mathbb P}_{\infty} = \int_{\Omega} d{\bf r}\ u_{st.}({\bf r}),
\end{equation}
which means that there is a finite probability that the particle will never 
be absorbed at the boundaries, resulting in a divergent MFPT.
This is a general result valid whenever the $g$-function has a finite asymptotic limit, 
corresponding to a slowing dynamics ending in a ``frozen'' particle distribution.


\section{Conclusions}

In this paper we have considered the $g$-fractional diffusion in $\mathbb{R}^d$ 
with absorbing (Dirichlet) boundary conditions. 
This generalizes the previous one-dimensional study \cite{PRE2023} to the $d$-dimensional case.
We show that it is possible to find the explicit representation of the solution for a generic $g$-function and bounded domain $\Omega$ in $\mathbb{R}^d$. 
 This is the first general treatment of $g$-fractional diffusion
in generic $d$ dimension and we obtain a general representation that can be
used in realistic models, beyond the one-dimensional analysis developed
in the previous literature.
We have then analyzed the FPTD and its dependence on 
the particular choice of the function $g$, leading to a nontrivial generalization of the classical Caputo-fractional diffusion in a bounded domain. We have considered the interesting case of the rectangular domain, obtaining the exact form of the solution.
Finally, we have devoted a section to a physical application 
related to the Dodson diffusion equation recently considered in \cite{garra}.
In previous research on this topic, the authors obtained the fundamental solution in the one-dimensional case, while here we derive the solution of the Dirichlet problem in higher dimensions.
This allows us to discuss the effects of choosing a bounded $g$-function,
resulting in a finite asymptotic survival probability and undefined MFPT.

In conclusion, we have demonstrated the nontrivial role of the $g$-function on diffusive behavior 
in $d$-dimensional domains and
its influence on the shape of FPTDs and the existence of finite MFPTs.
It would be interesting to extend the analysis to domains of different shapes 
(such as spherical or cylindrical)
and consider different boundary conditions, such as, for example, partial absorption 
\cite{AG_2020,ANG_2015}
with, possibly, time-dependent rates (see, for example, \cite{Bre_2022} and references therein)
, in order to study the combined effects of fractional diffusion and boundary properties on first-passage processes.


\section*{Acknowledgments}
L.A. acknowledge financial support from the Italian Ministry of University and Research (MUR) under the PRIN2020 Grant No. 2020PFCXPE.

    \end{document}